\documentclass{amsart}
\usepackage{amssymb,euscript,amsmath, mathrsfs}
\usepackage[dvips]{graphicx}
\usepackage[dvips]{color}
\usepackage{epsfig}

\newcounter{ENUM}
\newcommand{\itm}{\item}
\newenvironment{ilist}{\renewcommand{\theENUM}{\roman{ENUM}}\renewcommand{\itm}{\addtocounter{ENUM}{1}\item[(\theENUM)]}\begin{itemize}\setcounter{ENUM}{0}}{\end{itemize}}
\newenvironment{alist}[1][0]{\renewcommand{\theENUM}{\alph{ENUM}}\renewcommand{\itm}{\addtocounter{ENUM}{1}\item[\theENUM)]}\begin{itemize}\setcounter{ENUM}{#1}}{\end{itemize}}

\def\Z{{\mathbb Z}}
\def\N{{\mathbb N}}

\def\R{{\mathbb R}}

\def\cH{{\mathcal H}}

\def\L{{\mathcal L}}

\def\sgn{\operatorname{sign}}

\def\conv{\mathrm{conv}}
\def\vol{\mathrm{Vol}}

\newtheorem{thm}{Theorem}[section]
\newtheorem{prop}[thm]{Proposition}
\newtheorem{lem}[thm]{Lemma}
\newtheorem{cor}[thm]{Corollary}

\theoremstyle{definition}
\newtheorem{defn}[thm]{Definition}
\newtheorem{ques}[thm]{Question}
\newtheorem{ex}[thm]{Example}

\theoremstyle{remark}

\newtheorem{rem}[thm]{Remark}

\numberwithin{equation}{section}

\begin{document}
\title{Ehrhart Polynomials of Cyclic Polytopes}
\author{Fu Liu}
\begin{abstract}
The Ehrhart polynomial of an integral convex polytope counts the number of lattice points in dilates of the polytope. In \cite{conj}, the authors conjectured that for any cyclic polytope with integral parameters, the Ehrhart polynomial of it is equal to its volume plus the Ehrhart polynomial of its lower envelope and proved the case when the dimension $d = 2.$ In our article, we prove the conjecture for any dimension. 
\end{abstract}
\maketitle

\section{Introduction}
For any \textit{integral convex polytope} $P,$ that is to say, a convex polytope whose vertices have integral coordinates, any positive integer $m \in \N,$ we use $i(P, m)$ to denote the number of lattice points in $mP,$ where $mP = \{ mx \ | \ x \in P \}$ is a \textit{dilated polytope} of $P.$ In our paper, we will focus on a special class of polytopes, cyclic polytopes, which are defined in terms of the moment curve:

\begin{defn}
The {\bf moment curve} in $\R^d$ is defined by $$\nu_d: \R \to \R^d, t \mapsto \nu_d(t) =  \left( \begin{array}{cc} t \\ t^2 \\ \vdots \\ t^d \end{array} \right).$$
Let $T = \{t_1, \dots, t_n \}_<$ be a linear ordered set. Then the {\bf cyclic polytope} $C_d(T) = C_d(t_1, \dots, t_n)$ is the convex hull $\conv \{v_d(t_1), v_d(t_2), \dots, v_d(t_n) \}$ of $n > d$ distinct points $\nu_d(t_i), 1 \le i \le n,$ on the moment curve.
\end{defn}

The main theorem in our article is the one conjecured in \cite[Conjecture 1.5]{conj}: 

\begin{thm}\label{main}
For any integral cyclic polytope $C_d(T),$ 
$$i(C_d(T), m) = \vol(m C_d(T)) + i(C_{d-1}(T), m).$$ Hence,
$$i(C_d(T), m) = \sum_{k=0}^d \vol_k(m C_k(T)) = \sum_{k=0}^d \vol_k(C_k(T)) m^k,$$ where $\vol_k(m C_k(T))$ is the volume of $m C_k(T)$ in $k$-dimensional space, and we let $\vol_0(m C_0(T)) = 1.$
\end{thm}

One direct result of Theorem \ref{main} is that $i(C_d(T), m)$ is always a polynomial in $m.$ This result was already shown by Eug\`{e}ne Ehrhart for any integral polytope in 1967. Thus, we call $i(P,m)$ the Ehrhart polynomial of $P$ when $P$ is an integral polytope. There is much work on the coefficients of Ehrhart polynomials. For instance it's well known that the leading and second coefficients of $i(P,m)$ are the normalized volume of $P$ and one half of the normalized volume of the boundary of $P.$ But there is no known explicit method of describing all the coefficients of Ehrhart polynomials of general integral polytopes. However, because of some special properties that cyclic polytopes have, we are able to calculate the Ehrhart polynomial of cyclic polytopes in the way described in Theorem \ref{main}.

In this paper, we use a standard triangulation decomposition of cyclic polytopes, and careful counting of lattice points to reduce Theorem \ref{main} to the case n=d+1, (Theorem \ref{main2}).
We then prove Theorem \ref{main2} with the use of certain linear transformations and decompositions of polytopes containing our cyclic polytopes.

\section{Statements and Proofs}
All polytopes we will consider are full-dimensional, so for any convex 
polytope $P,$ we use $d$ to denote both the dimension of the ambient 
space $\R^d$ and the dimension of $P.$ Also, We use $\partial P$ and 
$I(P)$ to denote the boundary and the interior of $P,$ respectively.

For simplicity, for any region $R \subset \R^d,$ we denote by $\L(R) := R \cap \Z^d$ the set of lattice points in $R.$

Consider the projection $\pi: \R^d \to \R^{d-1}$ that forgets the last coordinate. In \cite[Lemma 5.1]{conj}, the authors showed that the inverse image under $\pi$ of a lattice point $y \in C_{d-1}(T) \cap \Z^{d-1}$ is a line that intersects the boundary of $C_d(T)$ at integral points. By a similar argument, it's easy to see that this is true when we replace the cyclic polytopes by their dilated polytopes. Note that $\pi(m C_d(T)) = m C_{d-1}(T),$ so for any lattice point $y$ in $m C_{d-1}(T)$ the inverse image under $\pi$ intersects the boundary at lattice points.

\begin{defn}
For any $x$ in a real space, let $l(x)$ denote the last coordinate of $x.$

For any polytope $P \subset \R^d$ and any point $y \in \R^{d-1},$ let $\rho(y, P) = \pi^{-1}(y) \cap P$ be the intersection of $P$ with the inverse image of $y$ under $\pi.$ Let $p(y, P)$ and $n(y, P)$ be the point in $\rho(y,P)$ with the largest and smallest last coordinate, respectively. If $\rho(y,P)$ is the empty set, i.e., $y \not\in \pi(P),$ then let $p(y, P)$ and $n(y,P)$ be empty sets as well. Clearly, $p(y, P)$ and $n(y, P)$ are on the boundary of $P.$
Also, we let $\rho^+(y,P) = \rho(y,P)\setminus n(y,P),$ and for any $S \subset \R^{d-1},$ $\rho^+(S,P) = \cup_{y \in S} \rho^+(y,P).$ 

Define $PB(P) = \bigcup_{y \in \pi(P)} p(y,P)$ to be the {\it positive boundary} of $P;$ $NB(P) = \cup_{y \in \pi(P)} n(y,P)$ to be the {\it negative boundary} of $P$ and $\Omega(P) = P \setminus NB(P) = \rho^+(\pi(P), P) = \cup_{y \in \pi(P)} \rho^+(y,P)$ to be the {\it nonnegative part} of $P.$

For any facet $F$ of $P,$ if $F$ has an interior point in the positive boundary of $P,$ (it's easy to see that $F \subset PB(P))$ then we call $F$ a {\it positive facet} of $P$ and define the sign of $F$ as $+1: \sgn(F) = + 1.$ Similarly, we can define the {\it negative facets} of $P$ with associated sign $-1.$
\end{defn}

It's not hard to see that $\pi$ induces a bijection of lattice points between $NB(m C_d(T))$ and  $\pi(m C_d(T)) = m C_{d-1}(T).$ Hence, Theorem \ref{main} is equivalent to the following Proposition:

\begin{prop}\label{main1}
$\vol(m C_d(T)) = |\L(\Omega(m C_d(T)))|.$
\end{prop}

From now on, we will consider any polytopes or sets as {\it multisets} which allow {\it negative multiplicities.} We can consider any element of a multiset as a pair $(x, m),$ where $m$ is the multiplicity of element $x.$ Then for any multisets $M_1, M_2$ and any integers $m,n$ and $i,$ we define the following operators:
\begin{ilist}
\itm Scalar product: $i M_1 = i \cdot M_1 = \{ (x, i m) \ | \ (x, m) \in M_1\}.$
\itm Addition: $M_1 \oplus M_2 = \{ (x, m + n) \ | \ (x, m) \in M_1, (x, n) \in M_2 \}.$
\itm Subtraction: $M_1 \ominus M_2 = M_1 \oplus ((-1) \cdot M_2).$
\end{ilist}

It's clear that the following holds:
\begin{lem} \
\begin{alist}
\itm $\forall R_1, \dots, R_k \subset \R^d,$ $\forall i_1,\dots, i_k \in \Z: \L(\bigoplus_{j=1}^k i_j R_j) = \bigoplus_{j=1}^k i_j \L(R_j). $
\itm For any polytope $P \subset \R^d,$  $\forall R_1, \dots, R_k \subset \R^{d-1},$ $\forall i_1,\dots, i_k \in \Z:$ $$\rho^+\left( \bigoplus_{j=1}^k i_j R_j, P \right) = \bigoplus_{j=1}^k i_j\rho^+( R_j, P). $$
\end{alist}
\end{lem}

Let $P$ be a convex polytope. For any $y$ an interior point of $\pi(P),$ since $\pi$ is a continous open map, the inverse image of $y$ contains an interior point of $P.$ Thus $\pi^{-1}(y)$ intersects the boundary of $P$ exactly twice; for any $y$ a boundary point of $\pi(P).$ Again because $\pi$ is an open map, we have that $\rho(y, P) \subset \partial P,$ so $\rho(y,P) = \pi^{-1}(y) \cap \partial P$ is either one point or a line segment. We hope that $\rho(y,P)$ always has only one point, so we define the following polytopes and discuss several properties of them. 

\begin{defn}
We call a convex polytope $P$ a \textit{nice} polytope with respect to $\pi$ if for any $y \in \partial \pi(P),$ $|\rho(y,P)| = 1$ and for any lattice point $y \in \pi(P),$ $\pi^{-1}(y)$ intersects $\partial P$ at lattice points. 
\end{defn}

\begin{lem}\label{pnp}
 A nice polytope $P$ has the following properties:
\begin{ilist}
\itm For any $y \in I(\pi(P)),$ $\pi^{-1}(y) \cap \partial P = \{p(y,P), n(y,P)\}.$ In particular, if $y$ is a lattice point, then $p(y,P)$ and $n(y,P)$ are each lattice points.
\itm For any $y \in \partial \pi(P),$ $\pi^{-1}(y) \cap \partial P = \rho(y,P) = p(y, P) = n(y, P),$ so $\rho^+(y,P) = \emptyset.$ In particular, when $y$ is a lattice point, $\rho(y,P)$ is a lattice point as well.
\itm $\L$ and $\rho^+$ commute: for any $R \subset \R^{d-1},$ $\L(\rho^+(R, P)) = \rho^+(\L(R), P).$ 
\itm Let $R$ be a region containing $I(\pi(P)).$ Then $$\Omega(P) = \rho^+(R, P) = \bigoplus_{y \in R} \rho^+(y, P).$$ Moreover, $$|\L(\Omega(P))| = \sum_{y \in \L(R)} l(p(y,P)) - l(n(y,P)).$$ (By convention, if $y \not\in \pi(P),$ we let $l(p(y,P)) - l(n(y,P)) = 0.$)
\itm If $P$ is decomposed into \it{nice} polytopes $P_1 \cup \cdots \cup P_k,$ then $\Omega(P) = \bigoplus_{i=1}^k \Omega(P_i),$ so $\L(\Omega(P)) = \bigoplus_{i=1}^k \L(\Omega(P_i)).$ 
\itm The set of facets of $P$ are partitioned into the set of positive facets and the set of negative facets, i.e., every facet is either positive or negative but not both.
\end{ilist}
\end{lem}

\begin{proof}
The first three and last properties are immediately true. And the fourth one follows directly from the second one. The fifth property can be checked by considering the definition of $\Omega.$ 
\end{proof}

By using these properties, we are able to give the following proposition about a nice convex polytope:

\begin{prop}\label{psum}
 Let $P$ be a nice convex polytope with respect to $\pi$ such $\pi(P)$ is also nice, and all the points in $P$ have nonnegative last coordinate. Suppose further that for any facet $F$ of $P,$ $\pi(F)$ is a nice polytope with respect to $\pi.$ Then
$$\Omega(P) = \bigoplus_{F: \mbox{ a facet of } P} \sgn(F) \rho^+(\Omega(\pi(F)), \conv(F, \pi(F))),$$
where $\conv(F, \pi(F))$ denotes the convex hull of the set $F \cup \{ (y', 0)' \ | \ y \in \pi(F) \},$ i.e. the region between $F$ and its projection onto the hyperplane $\{(x_1, \dots, x_d)' \ | \ x_d = 0 \}.$  (Note, for any vector $v,$ we use $v'$ to denote its transpose. So for a vertical vector $y,$ $(y',0)'$ is just the vector obtained from $y$ by attaching a zero to the bottom of $y.$)
\end{prop}

\begin{proof}
A special case of Lemma \ref{pnp}/(iv) is when $R = \Omega(\pi(P)),$ so we have 
$$\Omega(P) = \rho^+(\Omega(\pi(P)), P) = \bigoplus_{y \in \Omega(\pi(P))} \rho^+(y, P).$$

Now for any points $a$ and $b,$ we use $(a, b]$ to denote the half-open line segment between $a$(excluding) and $b$(including). Then, $\rho^+(y,P) = (n(y,P), p(y,P)] = ( ((y',0)', p(y, P)] \ominus ((y',0)', n(y,P)] ).$ Therefore,
\begin{align*}
\Omega(P) &=  \bigoplus_{y \in \Omega(\pi(P))} ( ((y',0)', p(y, P)] \ominus ((y',0)', n(y,P)] ) \\
&=(\bigoplus_{y \in \Omega(\pi(P))} ((y',0)', p(y, P)]) \bigoplus (\bigoplus_{y \in \Omega(\pi(P))} (-1)\cdot ((y',0)', n(y,P)] ).
\end{align*}

Let $F_1, F_2, \dots, F_\ell$ be all the positive facets of $P$ and $F_{\ell+1}, \dots, F_k$ be all the negative facets. Then it's clear that $\pi(F_1) \cup \pi(F_2) \cup \cdots \cup \pi(F_\ell)$ and $\pi(F_{\ell+1}) \cup \cdots \cup \pi(F_k)$ both give a decomposition of $\pi(P).$ Therefore by Lemma \ref{pnp}/(v), we have that $\Omega(\pi(P)) = \bigoplus_{i=1}^{\ell} \Omega(\pi(F_i)) = \bigoplus_{j=\ell+1}^k \Omega(\pi(F_j)).$ Hence,
\begin{eqnarray*}
\bigoplus_{y \in \Omega(\pi(P))} ((y',0)', p(y, P)] &=& \bigoplus_{i=1}^{\ell} \bigoplus_{y \in \Omega(\pi(F_i))} ((y',0)', p(y, P)] \\ 
&=& \bigoplus_{i=1}^{\ell} \rho^+(\Omega(\pi(F_i)), \conv(F_i, \pi(F_i))).
\end{eqnarray*}
Similarly, we will have
$$\bigoplus_{y \in \Omega(\pi(P))} (-1)\cdot ((y',0)', n(y,P)] = \bigoplus_{j=\ell+1}^k (-1) \rho^+(\Omega(\pi(F_j)), \conv(F_j, \pi(F_j))).$$
Thus, by putting them together, we get
$$\Omega(P) = \bigoplus_{\mbox{$F:$ a facet of $P$}} \sgn(F) \rho^+(\Omega(\pi(F)), \conv(F, \pi(F))).$$

\end{proof}

In the last proposition, we used a new notation $\conv(F, \pi(F))$ to denote certain polytopes. For polytopes that can be written in this way, we have the following lemma:

\begin{lem}\label{chp}
Let $H$ be a hyperplane in $\R^d,$ and let $S_1, S_2$ be two convex polytopes such that $S_1 \subset S_2 \subset H$ and the last coordinates of all of their points are nonnegative. Then for any $y \in \pi(S_1),$ $\rho^+(y, \conv(S_1, \pi(S_1))) = \rho^+(y, \conv(S_2, \pi(S_2))).$
\end{lem}

\begin{proof}
Trivial.
\end{proof}

After discussing some properties of nice polytopes with respect to $\pi,$ we come back to the dilated cyclic polytopes which are our main interest and show that they are nice:

\begin{lem} $m C_d(T)$ is a nice polytope with respect to $\pi.$ 
\end{lem}

\begin{proof}
We already argued that $m C_d(T)$ satisfies the second condition to be \textrm{nice.} So it's left to check that $|\rho(y, C_d(T))| = 1$ for any $y \in \partial C_{d-1}(T).$

Let $y = (y_1, y_2, \dots, y_{d-1})'$ and suppose $y$ is on a facet $F$ of $m C_{d-1}(T)$ and without loss of generality, let $m \nu_{d-1}(t_{1}), m \nu_{d-1}(t_{2}), \dots, m \nu_{d-1}(t_{{d-1}})$ be the $d-1$ vertices of $F.$ Then there exist $\lambda_1, \dots, \lambda_{d-1} \in \R_{\ge 0}$ such that $y = \sum_{j=1}^{d-1} \lambda_j m \nu_{d-1}(t_{j}) $ and $\sum_{j=1}^{d-1} \lambda_j = 1.$

Let $x \in \pi^{-1}(y) \cap m C_d(T).$ There exist $\lambda_1', \dots, \lambda_n' \in \R_{\ge 0}$ such that $x = \sum_{j=1}^n \lambda_j' m \nu_{d}(t_j)$ and $\sum_{j=1}^n \lambda_j' = 1.$ Then $y = \pi(x) = \sum_{j=1}^n \lambda_j' m \nu_{d-1}(t_j).$ Since $y$ is on the facet $F,$ $\lambda_j' = 0$ unless $1 \le j \le d - 1.$ Thus $y = \sum_{j=1}^{d-1} \lambda_j' m \nu_{d-1}(t_j)$ and $\sum_{j=1}^{d-1} \lambda_j' = 1.$ Therefore $\lambda_j = \lambda_j', 1 \le j \le d-1.$ Hence $x = \sum_{j=1}^{d-1} \lambda_j m \nu_{d}(t_{j})$ is the only point in $\pi^{-1}(y) \cap m C_d(T).$
\end{proof}

We know that for any cyclic polytope $C_d(T)$ with $n = |T| > d + 1,$ we can decompose it into $n - d$ cyclic polytopes $P_1 \cup \dots \cup P_{n-d},$ which is a triangulation of $C_d(T)$ and where $P_i$'s are all defined by $(d+1)$-element integer sets. Therefore by the fourth property in Lemma \ref{pnp}, we have that $\L(\Omega(C_d(T))) = \bigcup_{i=1}^{n-d} \L(\Omega(P_i)).$ Thus $|\L(\Omega(P))| = \sum_{i=1}^k |\L(\Omega(P_i))|.$ Note that we also have $\vol(C_d(T)) = \sum_{i=1}^{n-d} \vol(P_i).$ We conclude that to prove Proposition \ref{main1}, it is enough to prove the following:

\begin{thm}\label{main2}
For any integer sets $T$ with $n = |T| = d+1,$ $\vol(m C_d(T)) = |\L(\Omega(m C_d(T)))|.$ 
\end{thm}

\begin{defn}
A map $\varphi: \R^d \to \R^d$ is {\it structure perserving} if it preserves volume and it commutes with the following operations: 
\begin{ilist}
\itm $\L:$ taking lattice points of a region $R \subset \R^d;$
\itm $\conv:$ taking the convex hull of a collection of points;
\itm $\Omega:$ taking the nonnegative part of a convex polytope;
\itm $PB:$ taking the positive boundary of a convex polytope;
\itm $NB:$ taking the negative boundary of a convex polytope.
\end{ilist}
\end{defn}

\begin{rem} Here $\varphi$ commuting with $\conv$ implies (or is equivalent to) that for any set of points $x_1, \dots, x_k \in \R^d,$ and for any $\lambda_1, \dots, \lambda_k \in \R^{\ge 0}$ with $\sum_{i=1}^k \lambda_i = 1,$ $T(\sum_{i=1}^k \lambda_i x_i) = \sum_{i=1}^k \lambda_i T(x_i).$ Therefore $\varphi$ is an affine transformation, which can be defined by a $d \times d$ matrix $A$ and a vector $u \in \R^d: T(x) = Ax + u.$ Moreover, $\varphi$ commuting with $PB$ and $NB$ implies that $\varphi$ preserves the positive facets and negative facets of a convex polytope.
\end{rem}

\begin{lem}\label{ltp}
Let $A$ be a $d \times d$ integral lower triangular matrix with $1$'s on its diagonal, and $u$ be an integral vector in $\R^d.$ Then $\varphi: x \mapsto Ax + u$ gives a map which is structure preserving, and so does $\varphi^{-1}.$ Therefore, $\varphi$ is a bijection from $\Z^d$ to $\Z^d.$ Hence, for any subset $S \in \R^d,$ $|\L(S)| = |\L(\varphi(S))|.$

Moreover, for any $y \in \R^{d-1},$ if we define $\tilde{\varphi}(y) = \tilde{A}y + \tilde{u},$ where $\tilde{A}$ is the right upper $(d-1) \times (d-1)$ matrix of $A$ and $\tilde{u} = \pi(u),$ then $\rho^+(\tilde{\varphi}(y), \varphi(P)) = \varphi(\rho^+(y, P)),$ for any polytope $P.$

\end{lem}

\begin{proof}
The determinant of $A$ is $1,$ hence $\varphi$ is volume preserving. It's easy to check that $\varphi$ commutes with $\L$ and $\conv.$ To show that $\varphi$ commutes with $\Omega, PB$ and $NB,$ it suffices to show that for all $x_1, x_2 \in \R^d$ with $\pi(x_1) = \pi(x_2)$ and $l(x_1) > l(x_2),$ then $\pi(\varphi(x_1)) = \pi(\varphi(x_2))$ and $l(\varphi(x_1)) > l(\varphi(x_2)).$ This is not hard to check using the fact that $A$ is a lower triangular matrix with $1$'s on its diagonal. Hence, $\varphi$ is structure preserving.

Note that $\varphi^{-1}$ maps $x$ to $A^{-1}x - A^{-1}u.$ But we know that $A^{-1}$ is also an integral lower triangular matrix with $1$'s on its diagonal and $-A^{-1}u$ is an integral vector. So $\varphi^{-1}$ is structure preserving as well.

It's clear that $\tilde{\varphi} = \pi \circ \varphi \circ \pi^{-1},$ which implies that $\pi^{-1} \circ \tilde{\varphi} = \varphi \circ \pi^{-1}.$  So we have $$x \in \varphi(\rho^+(y,P)) \Leftrightarrow \varphi^{-1}(x) \in \rho^+(y,P) = \pi^{-1}(y) \cap P $$
$$\Leftrightarrow x \in \varphi(\pi^{-1}(y)) \cap \varphi(P) = \pi^{-1}(\tilde{\varphi}(y)) \cap \varphi(P) \Leftrightarrow x \in \rho^{+}(\tilde{\varphi}(y), \varphi(P)).$$
\end{proof}

Now for any real numbers $r_1, r_2, \dots, r_d,$ we consider the $d \times d$ lower triangular matrices 
\[
	A_{r_1,\dots, r_d}(i,j) =
	\left\{
	\begin{array}{ll}
	(-1)^{i-j} e_{i-j}(r_1, \dots, r_i), & i \ge j\\
	0, & i < j
	\end{array}
	\right.
	\] and
\[
	B_{r_1,\dots, r_d}(i,j) =
	\left\{
	\begin{array}{ll}
	1, & i = j \\
	0, & i \ne j \& i < d\\
	(-1)^{i-j} e_{i-j}(r_1, \dots, r_i), & j \ne i = d
	\end{array}
	\right.
	\]
 where $e_k(r_1, \dots, r_l) = \sum_{i_1 < i_2 < \cdots < i_k} r_{i_1} r_{i_2} \dots r_{i_k}$ is the $k$th elementary symmetric function in $r_1, \dots, r_l.$

For simplicity, we allow a map originally defined on $\R^d$ to work in higher dimension, by applying the map to the first $d$ coordinates. Then it's not hard to see that $A_{r_1,\dots, r_d} = A_{r_1,\dots, r_{d-1}} B_{r_1,\dots, r_d} = B_{r_1,\dots, r_d} A_{r_1,\dots, r_{d-1}}.$

We also define vectors 
$$u_{r_1, \dots, r_d} = \left(
	\begin{array}{ccccc}
	-r_1\\
	r_1 r_2\\
	-r_1 r_2 r_3\\
	\vdots \\
	(-1)^d r_1 r_2 \dots r_d
	\end{array}
	\right) = \left(
	\begin{array}{cc}
	-e_1(r_1)\\
	e_2(r_1,r_2)\\
	-e_3(r_1,r_2,r_3)\\
	\vdots \\
	(-1)^d e_d(r_1, r_2, \dots, r_d)
	\end{array}
	\right),$$
and
$$v_{r_1, \dots, r_d} = \left(
	\begin{array}{ccccc}
	0\\
	0\\
	\vdots \\
	0\\
	(-1)^d r_1 r_2 \dots r_d
	\end{array}
	\right) = \left(
	\begin{array}{cc}
	0\\
	0\\
	\vdots \\
	0\\
	(-1)^d e_d(r_1, r_2, \dots, r_d)
	\end{array}
	\right).
	$$
Similarly, we allow the addition operation between two vectors of different dimensions by adding the lower dimension one to the first corresponding coordinates of the higher one. Thus, $u_{r_1,\dots, r_d} = u_{r_1,\dots, r_{d-1}} +  v_{r_1,\dots, r_d}.$

Now we define maps $\varphi_{r_1, \dots, r_d}: x \mapsto A_{r_1, \dots, r_d} x + u_{r_1, \dots, r_d}$ and $ \phi_{r_1,\dots, r_d}: x \mapsto B_{r_1,\dots, r_d} x + v_{r_1,\dots, r_d}.$ Unlike $\varphi_{r_1, \dots, r_d},$ $\phi_{r_1, \dots, r_d}$ does not depend on the order of $r_i$'s. In other words, for any permutation $\sigma \in S_d,$ $\phi_{r_1, \dots, r_d} = \phi_{r_{\sigma(1)}, \dots, r_{\sigma(d)}}.$ 

Note that $\phi_{r_1, \dots, r_d}$ only changes the $d$th coordinate of a vector, so we have the following lemma:

\begin{lem}
$\varphi_{r_1, \dots, r_d} = \varphi_{r_1, \dots, r_{d-1}} \circ \phi_{r_1, \dots, r_d}.$ 
\end{lem}

\begin{rem}\label{phi}
When we consider $\varphi_{r_1, \dots, r_d}$ and $\phi_{r_1, \dots, r_d}$ operating on moment curve, we have
$$\varphi_{r_1, \dots, r_d}(\nu_d(t)) = A_{r_1, \dots, r_d}\left(
	\begin{array}{c}
	t \\
	t^2 \\
	\vdots \\
	t^d
	\end{array}
	\right) + u_{r_1, \dots, r_d} = \left(
	\begin{array}{c}
	(t - r_1)\\
	(t - r_1)(t - r_2)\\
	\vdots \\
	(t - r_1)(t - r_2) \cdots (t - r_d)\
	\end{array}
	\right),
$$
$$\phi_{r_1, \dots, r_d}(\nu_d(t)) = B_{r_1, \dots, r_d}\left(
	\begin{array}{c}
	t \\
	t^2 \\
	\vdots \\
	t^d
	\end{array}
	\right) + v_{r_1, \dots, r_d} = \left(
	\begin{array}{c}
	t\\
	t^2\\
	\vdots \\
	t^{d-1}\\
	(t - r_1)(t - r_2) \cdots (t - r_d)\
	\end{array}
	\right).
$$
\end{rem}

\begin{rem}
When $r_1, \dots, r_d$ are integers, $\varphi_{r_1, \dots, r_d}, \phi_{r_1, \dots, r_d}$ and their inverse maps are structure preserving by Lemma \ref{ltp}. 
\end{rem}

Now by using $\phi$'s (or $\varphi$'s), we are able to determine the sign of the facets of dilated cyclic polytopes:

\begin{prop}\label{cpf}
Let $P = m C_d(T),$ where $m \in \N$ and $T = \{t_1, t_2, \dots, t_n\}_<$ an integral ordered set. Let $F$ be a facet of $P$ determined by vertices $\nu_d(t_{i_1}), \nu_d(t_{i_2}), \dots, \nu_d(t_{i_d}).$ Let $k$ be the smallest element of the set $\{1, 2, \dots, n\} \setminus \{i_1, \dots, i_d\},$ then $\sgn(F) = (-1)^{d-k}.$ In particular, when $|T| = n = d+1,$ let $F_k$ be the facet of $P$ determined by all the vertices of $P$ except $\nu_d(t_{i_k}),$ then for $k \in [d], \sgn(F_k) = \sgn(\sigma_k),$ where $\sigma_k = (k,k+1,\cdots, d) \in S_d$ and $\sgn(F_{d+1}) = -1.$
\end{prop}

\begin{proof}
We first consider the case when $m = 1,$ i.e. $P$ is a cyclic polytope. Without loss of generality, we assume that $i_1 < i_2 < \cdots < i_d.$ Consider the polytope $Q = \phi_{t_{i_1}, \dots, t_{i_d}}(P).$ For $j = 1, 2, \dots, n,$ the last coordinate of the vertex of $Q$ which mapped from $\nu_d(t_j)$ is $l(\phi_{t_{i_1}, \dots, t_{i_d}}(\nu_d(t_j)))= (t_j - t_{i_1})(t_j - t_{i_2}) \cdots (t_j - t_{i_d}). $ Hence the last coordinates of the vertices of $\phi_{t_{i_1}, \dots, t_{i_d}}(F)$ are all $0$'s. So $\phi_{t_{i_1}, \dots, t_{i_d}}(F)$ is on the hyperplane obtained by setting the last coordinate to $0.$ Since $k$ is the smallest element not in $\{i_1, \dots, i_d \},$ $i_1 = 1, i_2 =2, \dots, i_{k-1} = k -1, i_k > k.$ So $t_k - t_{i_l} > 0$ when $l = 1, 2, \dots, k-1;$ and $t_k - t_{i_l} < 0$ when $l = k, k+1, \dots, d.$ Therefore $\sgn(l(\phi_{t_{i_1}, \dots, t_{i_d}}(\nu_d(t_k))) = (-1)^{d-k+1}.$ By using Gale's evenness condition \cite{Gale}, it's not hard to see that $\sgn(l(\phi_{t_{i_1}, \dots, t_{i_d}}(\nu_d(t_l))) = (-1)^{d-k+1},$ for all $ l \not\in \{i_1, \dots, i_d\}.$ Thus we can conclude that $l(\phi_{t_{i_1}, \dots, t_{i_d}}(P))$ is nonnegative if $d - k$ is odd, and is nonpositive if $d-k$ is even. Hence $\phi_{t_{i_1}, \dots, t_{i_d}}(F)$ and $F$ are negative facets if $d-k$ is odd, and positive facets if $d-k$ is even. So $\sgn(F) = (-1)^{d-k}.$ For $n = d+1,$ it's easy to see that $\sgn(\sigma_k) = (-1)^{d-k} = \sgn(F_k).$

For $m > 1,$ we just need consider the map $x \mapsto B_{t_{i_1}, \dots, t_{i_d}}x + m v_{t_{i_1}, \dots, t_{i_d}}$ instead of $\phi_{t_{i_1}, \dots, t_{i_d}},$ and then we will have similar results.

\end{proof}

\begin{lem}\label{lpr} 
For all $ d \in \R^+,$ for all $s_1, \dots, s_d \in \N,$ let $x_0 = 1$ and $P_{s_1, \dots, s_d} = \{(x_1, \dots, x_d) \in \R^d \ | \ \forall i \in [d]: 0 \le x_i \le s_i x_{i-1} \},$ $R_{s_1, \dots, s_d} = \Omega(P_{s_1, \dots, s_d}).$ Then $R_{s_1, \dots, s_d} = P_{s_1, \dots, s_d} \cap \{ x_d > 0 \}$ and for all $ d \ge 2: R_{s_1, \dots, s_d} = \rho^+(R_{s_1, \dots, s_{d-1}}, P_{s_1, \dots, s_d}).$

Moreover, the vertices of $P_{s_1, \dots, s_d}$ are 
$$\left(
	\begin{array}{c}
	0 \\
	0 \\
	\vdots \\
	0\\
	\end{array}
	\right),
\left(
	\begin{array}{c}
	s_1 \\
	0 \\
	\vdots \\
	0\\
	\end{array}
	\right),
\left(
	\begin{array}{c}
	s_1 \\
	s_1 s_2 \\
	0 \\
	\vdots \\
	0\\
	\end{array}
	\right), \dots,
\left(
	\begin{array}{c}
	s_1 \\
	s_1 s_2 \\
	s_1 s_2 s_3 \\
	\vdots \\
	s_1 s_2 \cdots s_d\\
	\end{array}
	\right)$$
and the positive boundary of $P_{s_1, \dots, s_d}$ is just the convex hull of the first $d-1$ vertices and the last one. Note the first $d-1$ vertices span a $(d-2)$-dimensional space $\{(x_1, \dots, x_d)' \ | \ x_d = x_{d-1} = 0 \}.$ Hence $PB(P_{s_1, \dots, s_d})$ is in the hyperplane spanned by this $(d-2)$-dimensional space and the last vertex.
\end{lem}

\begin{proof}
The first result is immediate by considering the definition of $\Omega.$

We have $R_{s_1, \dots, s_{d-1}} \subset P_{s_1, \dots, s_{d-1}},$ so 
\begin{align*}
\rho^+(R_{s_1, \dots, s_{d-1}}, P_{s_1, \dots, s_d}) \subset \rho^+(P_{s_1, \dots, s_{d-1}}, P_{s_1, \dots, s_d}) &= \rho^+(\pi(P_{s_1, \dots, s_{d}}), P_{s_1, \dots, s_d}) \\ &= \Omega(P_{s_1, \dots, s_d}) = R_{s_1, \dots, s_d}.
\end{align*}
 But for $\forall x =(x_1, \dots, x_d) \in R_{s_1, \dots, s_d},$ we have that $ x_d > 0 $ which implies that $s_d x_{d-1} > 0,$ so $x_{d-1} > 0.$ Therefore $\pi(x) \in R_{s_1, \dots, s_{d-1}}.$ Thus, $x \in \rho^+(R_{s_1, \dots, s_{d-1}}, P_{s_1, \dots, s_d}).$ Now we can conclude that $R_{s_1, \dots, s_d} = \rho^+(R_{s_1, \dots, s_{d-1}}, P_{s_1, \dots, s_d}).$

\end{proof}

\begin{thm}\label{pmain} Let $d \in \N$ and $T = \{ t_1, t_2, \dots, 
t_{d+1} \}_<$ be an integral ordered set, then $$\Omega(C_d(T)) = 
\bigoplus_{\sigma \in S_d} \sgn(\sigma) \varphi_{t_{\sigma(1)},\dots, 
t_{\sigma(d)}}^{-1}(R_{t_{d+1} - t_{\sigma(1)},\dots, t_{d+1} 
-t_{\sigma(d)}}).$$

\end{thm}

\begin{proof} We proceed by induction on $d.$
When $d = 1,$ $C_d(T)$ is just the interval $[t_1, t_2].$ Then the only element $\sigma \in S_1$ is the identity map. $R_{t_2 - t_1} = (0, t_2 - t_1].$ And $\varphi_{t_1}: x \mapsto x - t_1,$ so $\varphi_{t_1}^{-1}: x \mapsto x + t_1.$ Thus $\varphi_{t_1}^{-1}((0, t_2-t_1]) = (t_1, t_2] = \Omega([t_1, t_2]).$

Now we assume the theorem is true for dimensions less than $d,$ and we will prove the case of dimension $d (\ge 2).$ Let $P = \phi_{t_1, \dots, t_d} (C_d(T)),$ and let $v_i = \phi_{t_1, \dots, t_d}(\nu_d(t_i)), i \in [d+1],$ be the vertices of $P.$ Then for $i \in [d],$ 
$v_i = \left(
	\begin{array}{c}
	\nu_{d-1}(t_i) \\
	0
	\end{array}
	\right)$ and for $i = d + 1,$ 
$v_{d+1} = \left(
	\begin{array}{c}
	\nu_{d-1}(t_{d+1}) \\
	\prod_{i=1}^d (t_{d+1} - t_i))
	\end{array}
	\right).$
Since $\prod_{i=1}^d (t_{d+1} - t_i))>0,$ the last coordinates of all the points in $P$ are nonnegative. By Proposition \ref{psum}, we have that 
$$\Omega(P) = \bigoplus_{F: \mbox{ a facet of } P} \sgn(F) \rho^+(\Omega(\pi(F)), \conv(F, \pi(F))).$$

As in Proposition \ref{cpf}, we let $F_k$ be the facet of $C_d(T)$ determined by all the vertices of $C_d(T)$ except $\nu_d(t_{i_k}),$ 
then
$$\Omega(P) = \bigoplus_{k \in [d+1]} \sgn(\phi_{t_1, \dots, t_d}(F_k)) \rho^+(\Omega(\pi(\phi_{t_1, \dots, t_d}(F_k))), \conv(\phi_{t_1, \dots, t_d}(F_k), \pi(\phi_{t_1, \dots, t_d}(F_k)))).$$

For $k = d+ 1,$ $\tilde{F} = \phi_{t_1, \dots, t_d}(F_{d+1}) = \conv(\{v_i\}_{i=1}^d)$ is on the hyperplane $H_0 = \{(x_1,\dots, x_d)' \in \R^d \ | \ x_d = 0 \}.$ So $\conv(\tilde{F}, \pi(\tilde{F}))$ is just $\tilde{F}.$ Thus $\rho^+(\Omega(\pi(\tilde{F})), \conv(\tilde{F}, \pi(\tilde{F})))$ is an empty set.

And for $k \in [d],$ by Proposition \ref{cpf}, $\sgn(F_k) = \sgn(\sigma_k),$ where $\sigma_k = (k,k+1,\cdots, d) \in S_d.$ Let $T_k = T \setminus \{t_k\},$ then $\pi(\phi_{t_1, \dots, t_d}(F_k)) = \pi(F_k) = C_{d-1}(T_k),$ because $\phi_{t_1, \dots, t_d}$ just changes the last coordinates. It's easy to see that 
$$\conv(\phi_{t_1, \dots, t_d}(F_k), \pi(\phi_{t_1, \dots, t_d}(F_k))) = \conv(\{ v_i \}_{i\ne k} \cup \{v_{d+1}'\}),$$
 where $v_{d+1}' = \left(
	\begin{array}{c}
	\nu_{d-1}(t_{d+1}) \\
	0
	\end{array}
	\right)$ is the projection of $v_{d+1}$ to the hyperplane $H_0.$

Hence,
$$\Omega(P) = \bigoplus_{k \in [d]} \sgn(\sigma_k) \rho^+(\Omega(C_{d-1}(T_k)), \conv(\{ v_i \}_{i\ne k} \cup \{v_{d+1}'\})).$$

For any $k \in [d],$ $T_k = \{t_{\sigma_k(1)}, t_{\sigma_k(2)}, \dots, t_{\sigma_k(d-1)}, t_{d+1} \}_<$. By the induction hypothesis, we have that 
$$\Omega(C_{d-1}(T_k)) = \bigoplus_{\tau \in S_{d-1}} \sgn(\tau) \varphi_{t_{\sigma_k(\tau(1))},\dots, t_{\sigma_k(\tau(d-1))}}^{-1}(R_{t_{d+1} - t_{\sigma_k(\tau(1))},\dots, t_{d+1} -t_{\sigma_k(\tau(d-1))}}).$$
So,
\begin{eqnarray*}
& &\sgn(\sigma_k) \phi_{t_1,\dots, t_d}^{-1} (\Omega(C_{d-1}(T_k))) \\
&=& \bigoplus_{\tau \in S_{d-1}} \sgn(\sigma_k) \sgn(\tau) \varphi_{t_{\sigma_k(\tau(1))},\dots, t_{\sigma_k(\tau(d-1))}}^{-1}(R_{t_{d+1} - t_{\sigma_k(\tau(1))},\dots, t_{d+1} -t_{\sigma_k(\tau(d-1))}}) \\
&{=}& \bigoplus_{\sigma \in S_{d}: \sigma(d) =k} \sgn(\sigma) \varphi_{t_{\sigma(1)},\dots, t_{\sigma(d-1)}}^{-1}(R_{t_{d+1} - t_{\sigma(1)},\dots, t_{d+1} -t_{\sigma(d-1)}}). \quad (\mbox{let } {\sigma = \sigma_k \tau})
\end{eqnarray*}

Let $H_k$ be the hyperplane determined by $\phi_{t_1, \dots, t_d}(F_k),$ and $H_k^+ = \{ x \in H_k \ | \ l(x) \ge 0 \}.$
We claim that for all $\sigma \in S_d$ with $\sigma(d) = k,$ we have
$$ \varphi_{t_{\sigma(1)},\dots, t_{\sigma(d-1)}}^{-1}(PB(P_{t_{d+1} - t_{\sigma(1)},\dots, t_{d+1} -t_{\sigma(d-1)}, t_{d+1}-t_{\sigma(d)}})) \subset H_k^+.$$

Given this, we can pick a convex polytope $S_k \subset H_k,$ such that 
\begin{alist}
\itm The last coordinates of the points in $S_k$ are nonnegative;
\itm  $S_k$ contains $\varphi_{t_{\sigma(1)},\dots, t_{\sigma(d-1)}}^{-1}(PB(P_{t_{d+1} - t_{\sigma(1)},\dots, t_{d+1} -t_{\sigma(d-1)}, t_{d+1}-t_{\sigma(d)}})),$  for all $\sigma \in S_d$ with $\sigma(d) = k;$
\itm $S_k$ contains $\phi_{t_1, \dots, t_d}(F_k).$ 
\end{alist}
Hence, by Lemma \ref{chp}
\begin{eqnarray*}
& &\Omega(C_d(T)) \\
&=& \phi_{t_1, \dots, t_d}^{-1}(\Omega(P)) \\
&=& \bigoplus_{k \in [d]} \sgn(\sigma_k) \phi_{t_1, \dots, t_d}^{-1} (\rho^+(\Omega(C_{d-1}(T_k)), \conv(\{ v_i \}_{i\ne k} \cup \{v_{d+1}'\}))) \\
&=& \bigoplus_{k \in [d]} \sgn(\sigma_k) \phi_{t_1, \dots, t_d}^{-1} (\rho^+(\Omega(C_{d-1}(T_k)), \conv(S_k, \pi(S_k)))) \\
&=& \bigoplus_{k \in [d]} \sgn(\sigma_k) \phi_{t_1, \dots, t_d}^{-1} (\rho^+(\bigoplus_{\tau \in S_{d-1}} \sgn(\tau) \varphi_{t_{\sigma_k(\tau(1))},\dots, t_{\sigma_k(\tau(d-1))}}^{-1}(R_{t_{d+1} - t_{\sigma_k(\tau(1))},\dots, t_{d+1} -t_{\sigma_k(\tau(d-1))}}),\\
& & \qquad \qquad \qquad \qquad \qquad \indent \conv(S_k, \pi(S_k)))) \\
&=& \bigoplus_{k \in [d]} \bigoplus_{\sigma \in S_{d}, \sigma(d) = k} \sgn(\sigma) \phi_{t_1, \dots, t_d}^{-1} (\rho^+(  \varphi_{t_{\sigma(1)},\dots, t_{\sigma(d-1)}}^{-1}(R_{t_{d+1} - t_{\sigma(1)},\dots, t_{d+1} -t_{\sigma(d-1)}}), \\
& &\qquad \qquad \qquad \qquad \qquad \qquad \qquad \indent \varphi_{t_{\sigma(1)},\dots, t_{\sigma(d-1)}}^{-1}(P_{t_{d+1} - t_{\sigma(1)},\dots, t_{d+1} -t_{\sigma(d-1)}, t_{d+1}-t_{\sigma(d)}}))) \\
&=&  \bigoplus_{\sigma \in S_{d}}  \sgn(\sigma) \phi_{t_1, \dots, t_d}^{-1} \varphi_{t_{\sigma(1)},\dots, t_{\sigma(d-1)}}^{-1}(\rho^+( R_{t_{d+1} - t_{\sigma(1)},\dots, t_{d+1} -t_{\sigma(d-1)}}, \\
& &\qquad \qquad \qquad \qquad \qquad \qquad \qquad \indent \indent P_{t_{d+1} - t_{\sigma(1)},\dots, t_{d+1} -t_{\sigma(d-1)}, t_{d+1}-t_{\sigma(d)}})) \\
&=& \bigoplus_{\sigma \in S_{d}} \sgn(\sigma) \varphi_{t_{\sigma(1)},\dots, t_{\sigma(d)}}^{-1}(R_{t_{d+1} - t_{\sigma(1)},\dots, t_{d+1} -t_{\sigma(d-1)}, t_{d+1}-t_{\sigma(d)}}).
\end{eqnarray*}
Thus the claim implies the theorem.

Showing the claim is equivalent to showing that $$ PB(P_{t_{d+1} - t_{\sigma(1)},\dots, t_{d+1} -t_{\sigma(d-1)}, t_{d+1}-t_{\sigma(d)}}) \subset \varphi_{t_{\sigma(1)},\dots, t_{\sigma(d-1)}}(H_k^+).$$
Both $\varphi_{t_{\sigma(1)},\dots, t_{\sigma(d-1)}}$ and its inverse only work on the first $d-1$ coordinates of any point in $\R^d.$ Thus $\varphi_{t_{\sigma(1)},\dots, t_{\sigma(d-1)}}(H_k^+)$ is just $\varphi_{t_{\sigma(1)},\dots, t_{\sigma(d-1)}}(H_k) \cap \{x \in \R^d \ | \ l(x) \ge 0 \}.$ But it's clear that $PB(P_{t_{d+1} - t_{\sigma(1)},\dots, t_{d+1} -t_{\sigma(d-1)}, t_{d+1}-t_{\sigma(d)}})$ is in $\{x \in \R^d \ | \ l(x) \ge 0 \}.$ So it's enough to show that $$PB(P_{t_{d+1} - t_{\sigma(1)},\dots, t_{d+1} -t_{\sigma(d-1)}, t_{d+1}-t_{\sigma(d)}}) \subset \varphi_{t_{\sigma(1)},\dots, t_{\sigma(d-1)}}(H_k).$$

By Lemma \ref{lpr},  $PB(P_{t_{d+1} - t_{\sigma(1)},\dots, t_{d+1} -t_{\sigma(d-1)}, t_{d+1}-t_{\sigma(d)}})$ lies in the hyperplane $H$ which is spanned by $\{(x_1, \dots, x_d)' \ | \ x_d = x_{d-1} = 0 \}$ and 
$\left(
	\begin{array}{c}
	t_{d+1} - t_{\sigma(1)} \\
	(t_{d+1} - t_{\sigma(1)}) (t_{d+1} - t_{\sigma(2)}) \\
	(t_{d+1} - t_{\sigma(1)}) (t_{d+1} - t_{\sigma(2)}) (t_{d+1} - t_{\sigma(3)}) \\
	\vdots \\
	(t_{d+1} - t_{\sigma(1)}) (t_{d+1} - t_{\sigma(2)}) \cdots (t_{d+1} - t_{\sigma(d)})\\
	\end{array}
	\right)
.$ 
So we need show that $\varphi_{t_{\sigma(1)},\dots, t_{\sigma(d-1)}}(H_k) = H.$ Since $H_k$ is the hyperplane containing $\phi_{t_1, \dots, t_d}(F_k),$ it's enough to show that $\varphi_{t_{\sigma(1)},\dots, t_{\sigma(d-1)}}(\phi_{t_1, \dots, t_d}(F_k)) = \varphi_{t_{\sigma(1)},\dots, t_{\sigma(d)}}(F_k)$ is contained in $H.$ However, $F_k = \conv(\nu_d(T_k)).$ Meanwhile, by remark \ref{phi}, we have
$$\varphi_{t_{\sigma(1)},\dots, t_{\sigma(d)}}(\nu_d(t)) = \left(
	\begin{array}{c}
	(t - t_{\sigma(1)})\\
	(t - t_{\sigma(1)})(t - t_{\sigma(2)})\\
	\vdots \\
	(t - t_{\sigma(1)})(t - t_{\sigma(2)}) \cdots (t - t_{\sigma(d)})\
	\end{array}
	\right).
$$
Since $\sigma(d) = k,$ for any $i \in [d], i \ne k,$ $\varphi_{t_{\sigma(1)},\dots, t_{\sigma(d)}}(\nu_d(t_i))$ has the last two coordinates equal to $0.$ And for $i = d+1,$ $\varphi_{t_{\sigma(1)},\dots, t_{\sigma(d)}}(\nu_d(t_{d+1}))$ is exactly the last vertex of $P_{t_{d+1} - t_{\sigma(1)},\dots, t_{d+1} -t_{\sigma(d-1)}, t_{d+1}-t_{\sigma(d)}},$ which completes the proof the claim and hence the theorem.

\end{proof}

\begin{rem}
If we define $\varphi_{m, r_1, \dots, r_d}: x \mapsto A_{r_1, \dots, r_d} x + m u_{r_1, \dots, r_d},$ then similarly we can prove that 
$$\Omega(m C_d(T)) = \bigoplus_{\sigma \in S_d} \sgn(\sigma) \varphi_{m, t_{\sigma(1)},\dots, t_{\sigma(d)}}^{-1}(m R_{t_{d+1} - t_{\sigma(1)},\dots, t_{d+1} -t_{\sigma(d)}}).$$
\end{rem}

\begin{cor}\label{liden}
$$\L(\Omega(m C_d(T))) = \bigoplus_{\sigma \in S_d} \sgn(\sigma) \L(\varphi_{m,t_{\sigma(1)},\dots, t_{\sigma(d)}}^{-1}(m R_{t_{d+1} -t_{\sigma(1)},\dots, t_{d+1} -t_{\sigma(d)}})).$$ Hence,
$$|\L(\Omega(m C_d(T)))| = \sum_{\sigma \in S_d} \sgn(\sigma) |\L(m R_{t_{d+1} -t_{\sigma(1)},\dots, t_{d+1} -t_{\sigma(d)}})|.$$

\end{cor}

It's easy to see that $m R_{s_1, \dots, s_d} = R_{m s_1, s_2, \dots, s_d}.$ Moreover,
 $$|L(R_{s_1, \dots, s_d})| = \sum_{x_1 = 1}^{s_{1}} \sum_{x_2 = 1}^{s_{2} x_1} \hdots \sum_{x_n = 1}^{s_{n} x_{n-1}} 1.$$ Therefore, it's natural to look at the following:

\begin{lem}\label{lemh} For any nonnegative integers $a_1, a_2, \dots, a_n, $ let $$h(a_1, a_2, \dots, a_n) = \sum_{x_1 = 1}^{a_{1}} \sum_{x_2 = 1}^{a_{2} x_1} \hdots \sum_{x_n = 1}^{a_{n} x_{n-1}} 1.$$ Then  the only highest degree term of $h$ is ${\frac{1}{n!}} a_1^n a_2^{n-1} a_3^{n-2} \hdots a_n.$ This is also true when we consider $h$ as a polynomial just in the variable $a_1.$
\end{lem}

\begin{proof}[Proof of Lemma \ref{lemh}:] We will prove it by induction on $n.$

When $n = 1,$ $h(a_1) = \sum_{x_1 = 1}^{a_1}1 = a_1.$ Thus the lemma holds.

Assume the lemma is true for $n,$ and note that $h(a_1, a_2, \dots, a_{n+1}) = \sum_{x_1=1}^{a_1} h(a_2 x_1, a_3, \dots, a_{n+1}).$ By assumption, ${\frac{1}{n!}} a_2^n a_3^{n-1} \hdots a_{n+1} x_1^n$ is the only highest degree term of $h(a_2 x_1, a_3, \dots, a_{n+1})$ when we consider it as polynomial both in $y = a_2 x_1, a_3, \dots, a_{n+1}$ and in $y.$ This implies that ${ \frac{1}{n!}}a_2^n a_3^{n-1} \hdots a_{n+1} x_1^n$ is the only highest degree term of $h(a_2 x_1, a_3, \dots, a_{n+1})$ when we consider it both in $a_2, a_3, \dots, a_{n+1}$ and in $x_1.$ Then our lemma immediately follows from the fact that the highest degree term of $\sum_{x_1=1}^{a_1} x_1^n$ is ${\frac{1} {n+1}} a_1^{n+1}.$

\end{proof}

\begin{prop}\label{iden} For any nonnegative integers $a_1, a_2, \dots, a_n,$ let 
$\cH_m (a_1, a_2, \dots, a_n) = \sum_{\sigma \in S_n} \sgn(\sigma) h({m a_{\sigma(1)}}, {a_{\sigma(2)}} \dots, {a_{\sigma(n)}}).$ Then 
$$\cH_m(a_1, a_2, \dots, a_n) = {\frac{m^n} {n!}} \prod_{i=1}^n a_i \prod_{1 \le i < j \le n} (a_i - a_j).$$
\end{prop}

\begin{proof}[Proof of Proposition \ref{iden}:]\

Clearly if any of $a_i$'s is $0,$ then $\cH_m(a_1, \dots, a_n) = 0.$ Also for $1 \le i < j \le n,$ $\cH_m$ changes sign when we switch $a_i$ and $a_j$, i.e.,
$$\cH_m( \dots, a_i, \dots, a_j, \dots) = -\cH_m (\dots, a_j, \dots, a_i, \dots).$$
Therefore, $\cH_m(a_1, \dots, a_n)$ must be a multiple of 
$$\prod_{i=1}^n a_i \prod_{1 \le i < j \le n} (a_i - a_j),$$ which has degree $\frac{1}{2}n(n+1).$
 
So now it's enough to show that $\cH_m(a_1, \dots, a_n)$ is of degree $\frac{1}{2}n(n+1)$ and the coefficient of $a_1^n a_2^{n-1} a_3^{n-2} \hdots a_n$ in $\cH_m(a_1, \dots, a_n)$ is ${\frac{m^n}{n!}},$ which follows from Lemma \ref{lemh}.
\end{proof}

\begin{proof}[Proof of Theorem \ref{main2}:]\
By Corollary \ref{liden}, 
\begin{eqnarray*}
|\L(\Omega(m C_d(T)))| &=& \sum_{\sigma \in S_d} \sgn(\sigma) |\L(m R_{t_{d+1} -t_{\sigma(1)},\dots, t_{d+1} -t_{\sigma(d)}})| \\
&=& \cH_m(t_{d+1} -t_{\sigma(1)},t_{d+1}-t_{\sigma(2)}, \dots, t_{d+1} -t_{\sigma(d)}) \\
&=& {\frac{m^d}{d!}} \prod_{i=1}^d (t_{d+1} - t_i \prod_{1 \le i < j \le d} (t_i - t_j) \\
&=& {\frac{m^d}{d!}} \prod_{1 \le i < j \le d+1} (t_i - t_j) = \vol(mC_d(T)).
\end{eqnarray*}
\end{proof}

As we argued earlier in our paper, the proof of Theorem \ref{main2} completes the proof of Proposition \ref{main1} and thus proof of our main Theorem \ref{main}.

\section{Examples and Question}

In this section, we are going to show some examples to make some of the statements or their proofs in the last section more clear. We will use the cyclic polytope $P = C_d(T),$ where $d = 3, T = \{1,2,3,4\}$ throughout this section. Let $v_k = \nu_d(k)$ be the $i$th vertex of $P$ and $F_k = \conv(\{v_1,v_2,v_3,v_4\} \setminus v_k).$

\begin{ex}\label{excpf}
According to Proposition \ref{cpf}, $\sgn(F_1) = \sgn((1,2,3)) = 1, \sgn(F_2) = \sgn((2,3)) = -1, \sgn(F_3) = \sgn((3)) = 1$ and $\sgn(F_4) = -1.$ So $F_1$ and $F_3$ are positive facets while $F_2$ and $F_4$ are negative facets. $\Omega(P) = P \setminus (F_2 \cup F_4).$ 
\end{ex}

\begin{ex}[Example of Structure Preserving map] \

$\phi_{1,2,3}: x \mapsto \left(
	\begin{array}{ccc}
	1 & 0 & 0 \\
	0 & 1 & 0 \\
	11 & -6 & 1
	\end{array}
	\right) x + \left(
	\begin{array}{c}
	0 \\ 0 \\ -6 \end{array} \right),$ 
so $\phi_{1,2,3}(\nu_3(t) ) = \left( \begin{array}{c} t \\ t^2 \\ (t-1)(t-2)(t-3) \end{array} \right).$ In particular, $\phi_{1,2,3}(v_1) = (1,1,0)', \phi_{1,2,3}(v_2) = (2,4,0)',\phi_{1,2,3}(v_3) = (3,9,0)',$ and $\phi_{1,2,3}(v_4) = (4,16,6)'.$ Therefore, $\phi_{1,2,3}(P) = \conv((1,1,0)', (2,4,0)', (3,9,0)', (4,16,6)').$ Because $\phi_{1,2,3}$ is structure preserving, $\sgn(F_k') = \sgn(F_k),$ where $F_k' = \phi_{1,2,3}(F_k).$
\end{ex}

\begin{ex}\label{expsum}
$\phi_{1,2,3}(P)$ is a polytope that satisfies the hypothesis in the Proposition \ref{psum}, so we should have
$$\Omega(\phi_{1,2,3}(P)) = \bigoplus_{k \in [4]} \sgn(F_k') \rho^+(\Omega(\pi(F_k')), \conv(F_k', \pi(F_k'))).$$ Now we check it:

$\conv(F_4', \pi(F_4')) = \conv((1,1,0)',(2,4,0)',(3,9,0)').$ It's just a triangle in the hyperplane $H_0 = \{(x_1,x_2,x_3)' \ | \ x_3 = 0 \}.$ Therefore, $$\rho^+(\Omega(\pi(F_4')), \conv(F_4', \pi(F_4'))) = \emptyset.$$

$\pi(F_1') = \conv((2,4)',(3,9)',(4,16)') = C_2(2,3,4)$ is a triangle, whose positive facet is $\conv(\nu_2(2),\nu_2(4))$ and negative facets are $\conv(\nu_2(2), \nu_2(3))$ and $\conv(\nu_2(3), \nu_2(4)).$ Hence,$$\Omega(\pi(F_1')) = C_2(2,3,4) \setminus (\conv(\nu_2(2), \nu_2(3)) \cup \conv(\nu_2(3), \nu_2(4))),$$
$$\conv(F_1', \pi(F_1')) = \conv((2,4,0)',(3,9,0)',(4,16,0)',(4,16,6)').$$
Therefore,
\begin{eqnarray*}
& &\rho^+(\Omega(\pi(F_1')), \conv(F_1', \pi(F_1'))) \\
&=& \Omega(\conv((2,4,0)',(3,9,0)',(4,16,0)',(4,16,6)')) \ominus \Omega(\conv((3,9,0)',(4,16,0)',(4,16,6)')).
\end{eqnarray*}
Similarly, 
\begin{eqnarray*}
& &\rho^+(\Omega(\pi(F_2')), \conv(F_2', \pi(F_2'))) \\
&=& \Omega(\conv((1,2,0)',(3,9,0)',(4,16,0)',(4,16,6)')) \ominus \Omega(\conv((3,9,0)',(4,16,0)',(4,16,6)')).
\end{eqnarray*}
And 
\begin{eqnarray*}
& &\rho^+(\Omega(\pi(F_3')), \conv(F_3', \pi(F_3'))) \\
&=& \Omega(\conv((1,2,0)',(2,4,0)',(4,16,0)',(4,16,6)')) \ominus \Omega(\conv((2,4,0)',(4,16,0)',(4,16,6)')).
\end{eqnarray*}
Thus, 
\begin{eqnarray*}
& &\bigoplus_{k \in [4]} \sgn(F_k') \rho^+(\Omega(\pi(F_k')), \conv(F_k', \pi(F_k'))) = \bigoplus_{k \in [3]} \sgn(F_k') \rho^+(\Omega(\pi(F_k')), \conv(F_k', \pi(F_k'))) \\
&=&  \Omega(\conv((2,4,0)',(3,9,0)',(4,16,0)',(4,16,6)')) \ominus \Omega(\conv((3,9,0)',(4,16,0)',(4,16,6)')) \\
& & \ominus (\Omega(\conv((1,2,0)',(3,9,0)',(4,16,0)',(4,16,6)')) \ominus \Omega(\conv((3,9,0)',(4,16,0)',(4,16,6)'))) \\
& & \oplus (\Omega(\conv((1,2,0)',(2,4,0)',(4,16,0)',(4,16,6)')) \ominus \Omega(\conv((2,4,0)',(4,16,0)',(4,16,6)'))) \\
&=& \Omega(\conv((1,2,0)',(2,4,0)',(4,16,0)',(4,16,6)')) \ominus \Omega(\conv((2,4,0)',(4,16,0)',(4,16,6)')) \\
& & \oplus \Omega(\conv((2,4,0)',(3,9,0)',(4,16,0)',(4,16,6)')) \ominus \Omega(\conv((1,2,0)',(3,9,0)',(4,16,0)',(4,16,6)')) \\
&=& \Omega(\conv((1,2,0)',(2,4,0)',(3,9,0)',(4,16,0)',(4,16,6)')) \\
& & \ominus \Omega(\conv((1,2,0)',(3,9,0)',(4,16,0)',(4,16,6)')) \\
&=& \Omega(\phi_{1,2,3}(P))
\end{eqnarray*}
This agrees with Propostion \ref{psum}.

\end{ex}

We will illustrate explicitly how we get the formula in Theorem \ref{pmain} for $P = C_3(1,2,3,4)$ in the next example:
\begin{ex} 
Acoording to Propostion \ref{psum} or Example \ref{expsum}, we have:

$\Omega(\phi_{1,2,3}(P)) = \bigoplus_{k \in [3]} \sgn(F_k') \rho^+(\Omega(\pi(F_k')), \conv(F_k', \pi(F_k'))) $

$\Rightarrow \Omega(P) = \bigoplus_{k \in [3]} \sgn(F_k') \phi_{1,2,3}^{-1} \rho^+(\Omega(\pi(F_k')), \conv(F_k', \pi(F_k'))).$

$\phi_{1,2,3}^{-1} \rho^+(\Omega(\pi(F_1')), \conv(F_1', \pi(F_1'))) = \phi_{1,2,3}^{-1}(\Omega(\conv((2,4,0)',(3,9,0)',(4,16,0)',(4,16,6)')) \ominus \Omega(\conv((3,9,0)',(4,16,0)',(4,16,6)')) )
= \Omega(\conv((2,4,8)',(3,9,27)',(4,16,58)',(4,16,64)')) \ominus \Omega(\conv((3,9,27)',(4,16,58)',(4,16,64)')).$

$\varphi_{2,3,1}: x \mapsto \left(
	\begin{array}{ccc}
	1 & 0 & 0 \\
	-5 & 1 & 0 \\
	11 & -6 & 1
	\end{array}
	\right) x + \left(
	\begin{array}{c}
	-2 \\ 6 \\ -6 \end{array} \right),$ 
so $\varphi_{1,2,3}(\nu_3(t) ) = \left( \begin{array}{c} (t-2) \\ (t-2)(t-3) \\ (t-1)(t-2)(t-3) \end{array} \right).$
And $\varphi_{3,2,1}(x)$ is just $\varphi_{2,3,1}(x) - (1,0,0)'.$ Hence, $\varphi_{2,3,1}(\varphi_{3,2,1}^{-1}(x)) = x + (1,0,0)'.$ Therefore, 
\begin{eqnarray*}
& &\varphi_{2,3,1}(\phi_{1,2,3}^{-1} \rho^+(\Omega(\pi(F_1')), \conv(F_1', \pi(F_1')))) \\
&=& \Omega(\conv((0,0,0)',(1,0,0)',(2,2,0)',(2,2,6)')) \ominus \Omega(\conv((1,0,0)',(2,2,0)',(2,2,6)')) \\
&=& \Omega(\conv((0,0,0)',(2,0,0)',(2,2,0)',(2,2,6)')) \ominus \Omega(\conv((1,0,0)',(2,0,0)',(2,2,0)',(2,2,6)')) \\
&=& R_{2,1,3} \ominus \varphi_{2,3,1}(\varphi_{3,2,1}^{-1}(\Omega(\conv((0,0,0)',(1,0,0)',(1,2,0)',(1,2,6)')))) \\
&=& R_{2,1,3} \ominus \varphi_{2,3,1}(\varphi_{3,2,1}^{-1}(R_{1,2,3}))
\end{eqnarray*}
Thus, $$\phi_{1,2,3}^{-1} \rho^+(\Omega(\pi(F_1')), \conv(F_1', \pi(F_1'))) = \varphi_{2,3,1}^{-1}(R_{4-2,4-3,4-1}) \ominus \varphi_{3,2,1}^{-1}(R_{4-3,4-2,4-1}).$$
We will have similar results for $F_2$ and $F_3.$ Therefore, 
$$\Omega(P) = \bigoplus_{\sigma \in S_3} \sgn(\sigma) \varphi_{t_{\sigma(1)},t_{\sigma(2)}, t_{\sigma(3)}}^{-1}(R_{t_{4} - t_{\sigma(1)},t_{4}-t_{\sigma(2)}, t_{4} -t_{\sigma(3)}}),$$
which agrees to the Theorem \ref{pmain}.
\end{ex}

Now we will use  Theorem \ref{main} to calculate $i(C_3(1,2,3,4),m):$
\begin{ex}\label{ex}
According to Theorem \ref{main} $$i(C_3(T), m) = \sum_{k=0}^3 \vol_k(C_k(T)) m^k.$$
$C_3(T)$ itself is a simplex, so $$\vol_3(C_3(T)) = \frac{1}{3!} \prod_{1 \le i < j \le 4} (j-i) = 2.$$

$C_2(T)$ can be decomposed into simplices $C_2(1,2,3)$ and $C_2(1,3,4),$ thus $$\vol_2(C_2(T)) = \vol_2(C_2(1,2,3)) + \vol_2(C_2(1,3,4))  = \frac{1}{2!}[(2-1)(3-1)(3-2) + (3-1)(4-1)(4-3)] = 4.$$

$C_1(T)$ is just an interval $[1,4],$ so $\vol_1(C_1(T)) = 4 - 1 = 3.$

Therefore, $i(C_3(1,2,3,4), m) = 2m^3 + 4m^2 + 3m + 1.$

\end{ex}

Since our theorem gives a nice form of Ehrhart polynomials of cyclic polytopes, it's natural to ask the following:

\begin{ques}
Are there other integral polytopes which have the same form of Ehrhart polynomials as cyclic polytopes? In other words, what kind of integral $d$-polytopes $P$ are there whose Ehrhart polynomials will be in the form of the following?
$$i(P, m) =  \vol(mP) + i(\pi(P), m) = \sum_{k=0}^d \vol_k(\pi^{(k)}(P)) m^k,$$
where $\pi^{(k)}$ is the map which ignores the last $k$ coordinates of a point. 
\end{ques}

\subsection*{Acknowledgements} I would like to thank Richard Stanley for showing me the conjecture in \cite{conj}.

\bibliographystyle{hamsplain}
\bibliography{gen}

\end{document}